\documentclass[11pt,dvips]{article}
\usepackage{latexsym}
 \usepackage{amsmath,amsthm,amsfonts,amssymb} 
 \usepackage{url}%utile pour les preprints avec url

\usepackage{enumerate}% \begin{enumerate}[cas 1]
\usepackage{mdwlist} % enumerate*
\usepackage{graphicx}
 \usepackage[all]{xy}

\newtheorem{thm}{Theorem}%[section]
\newtheorem{thmbis}{Theorem}
\newtheorem*{thm*}{Theorem}
 
\newtheorem*{dfn*}{Definition}

\newtheorem*{cor*}{Corollary}

\newtheorem*{prop*}{Proposition} 
\newtheorem*{properties*}{Properties}

\newtheorem*{lem*}{Lemma} 
\newtheorem{lembis}[thmbis]{Lemma} 
 
\newtheorem*{claim*}{Claim} 
 
\newtheorem*{fact*}{Fact}

\newtheorem*{qst*}{Question}

\newtheorem*{pb*}{Problem}

\theoremstyle{remark}
\newtheorem*{rem*}{Remark}

\newtheorem*{example*}{Example}

\renewcommand{\phi}{\varphi} 
\newcommand{\m} {^{-1}}

 \newcommand {\Ra} {\Rightarrow}

\newcommand{\inc}{\subset}
\newcommand {\Z} {{\mathbb {Z}}}   
\newcommand {\Q} {{\mathbb {Q}}}  

\begin{document}

\title{Unsolvability of the isomorphism problem for [free abelian]-by-free groups}
\author{ Gilbert Levitt}

 \date{}
\maketitle

\begin{abstract}  The isomorphism problem for  [free abelian]-by-free groups is unsolvable.
\end{abstract}

\section{Introduction}

This note was inspired by   Enric Ventura's talk at the 2007 Dortmund conference. I am grateful    to the organizers of the conference,   to E.\  Ventura for encouraging me to write this note, and to C.\ Miller and E.\ Ventura for helpful comments.

In his talk, Ventura described work by Bogopolski-Martino-Ventura \cite{BMV} showing that the conjugacy problem is unsolvable in certain  $\Z^4$-by-free groups (groups mapping onto a free group $F_n$ with kernel isomorphic to $\Z^4$). 
Their proof relies on the fact that 
$GL(4,\Z)$ contains
$F_2\times F_2$, and $F_2\times F_2$ has a finitely generated subgroup with unsolvable membership problem by a   construction due to Mihailova \cite{Mih}.

We elaborate on these ideas to show that the isomorphism problem for $\Z^4$-by-free groups is unsolvable (after this note was first posted on arXiv, Bruno Zimmermann pointed out his paper \cite{Zi} which contains a similar result for $\Z^4$-by-surface groups; his argument also applies to $\Z^4$-by-free groups).

\begin{thm}   There is no algorithm to decide whether two  groups of the form $\Z^4\rtimes F_{15}$ are isomorphic.
\end{thm}

The isomorphism problem was known to be unsolvable  for free-by-free groups  \cite{Mil}. It is solved for  [free abelian]-by-[free abelian] groups \cite{Se}, and open for free-by-[free abelian] groups. In fact, the isomorphism problem for mapping tori $F_n\rtimes _\alpha\Z$ of automorphisms of free groups is open in general. It is solved when $n=2$ \cite{BMV2}, when  there are only   finitely many $\alpha$-periodic conjugacy classes (using \cite{DG} and \cite{GL}), and when $\alpha$ has finite order in $Out(F_n)$ because $F_n\rtimes _\alpha\Z$ is a generalized Baumslag-Solitar group and \cite{Fo} applies.

Bogopolski-Martino-Ventura use Mihailova's construction to obtain
a finitely generated subgroup $A\inc GL(4,\Z)$ which is orbit undecidable: there is no algorithm to decide whether two given elements $x,y\in \Z^4$ are in the same $A$-orbit. 
In the proof of Theorem 1, we use Mihailova's construction  to establish another undecidability result for subgroups of $GL(4,\Z)$:

\begin{thm} 
There is no algorithm which takes as an input two finite subsets of $GL(4,\Z)$ and decides whether the subgroups generated by these subsets are conjugate.
\end{thm}

C.\ Miller pointed out that this theorem may be deduced from the proof of Theorem 4.5 of \cite{Mil2}, which asserts that the isomorphism problem for finitely generated subgroups of $F_2\times F_2$ is unsolvable. 

 We do not know whether the number 4 in Theorems 1 and 2 is optimal. An algorithm for deciding whether two finitely generated subgroups of $GL(2,\Z)$ are conjugate may be found in \cite{Me}.

\section{Proofs}

We recall    Mihailova's construction. Let $H$ be a 2-generated finitely presented group with unsolvable word problem; there is such an $H$ with 12 relations, see \cite{BMV}. Let $\pi:F_2\to H$ be the projection, and let $L=\{(u,v)\in F_2\times F_2\mid \pi(u)=\pi(v)\}$. Mihailova's  key observations are that  $L$ has unsolvable membership problem in $F_2\times F_2$, and $L$ is finitely generated (see \cite{Mil}); note that $L$ is not finitely presented \cite{Gru}. If $H$ has 12 relations, then $L$ is 14-generated. 

\begin{lembis} 
If $a\in F_2\times F_2$ and $L\inc a\m La$, then $L=a\m La$.
\end{lembis}

\begin{proof}
Let $a=(x,y)$. Then $L\inc a\m La$ translates to 
$$\pi(u)=\pi(v)\implies \pi(xux\m)=\pi(yvy\m).
$$
Letting $u=v$ shows that $\pi(x\m y)$ is central in $H$ and therefore
$$
\pi(xux\m)=\pi(yvy\m)\implies \pi(u)=\pi(x\m y)\pi(v)(\pi(x\m y))\m\implies \pi(u)=\pi(v).
$$
\end{proof}

From now on we shall view $F_2\times F_2$ as a subgroup of $GL(4,\Z)$, by first embedding $F_2\times F_2$ as a finite index subgroup of $GL(2,\Z) \times GL(2,\Z)$, and then embedding $GL(2,\Z) \times GL(2,\Z)$ into $GL(4,\Z)$ as block-diagonal matrices.

Fix  elements $h_1,\dots, h_p$   of $GL(4,\Z)$ such that $\langle h_1,\dots, h_p\rangle=L$. Given $h\in F_2\times F_2\inc GL(4,\Z)$, let $G_h$ be the semi-direct product $\Z^4\rtimes F_{p+1}$ defined by the  following action of $F_{p+1}$ on $\Z^4$: the first $p$ generators $t_1,\dots, t_p$ act as $h_1,\dots, h_p$ respectively, the last generator $t$ acts as $h$. We write $G_1$ for the group obtained when $h$ is trivial.

\begin{lembis}
Given $h\in F_2\times F_2$, the following are equivalent:
\begin{enumerate}
\item $h\in L$;
\item $G_h$ is isomorphic to $G_1$;
\item $\langle L,h\rangle $ and $L$ are conjugate subgroups of $GL(4,\Z)$.
\end{enumerate}
\end{lembis}

Since membership in $L$ is undecidable, this lemma implies Theorems 1 and 2   (note that $p$  may be taken to be 14). Let us  now prove  it.

$1\Ra 2$ is easy. If $h\in L$, it may be expressed as a word $w$ in terms  of the elements $h_i^{\pm 1} $. One constructs an isomorphism from $G_h$ to $G_1$ which is the identity on $\Z^4$ and $t_1,\dots, t_p$, and maps $t$ to $w(t_i^{\pm 1})\, t$.

For $2\Ra 3$, note that $\Z^4$ is the unique maximal normal abelian subgroup of $G_h$. If $G_h$ is isomorphic    to $G_1$, the isomorphism   maps $\Z^4$ to $\Z^4$ and therefore induces a conjugacy between the actions of $G_h$ and $G_1$ on $\Z^4$. It follows that $\langle L,h\rangle $ and $\langle L,1\rangle =L$ are conjugate.

For $3\Ra 1$, assume that $\langle L,h\rangle =g\m Lg $ for some $g\in GL(4,\Z)$. We will show that $g$ has a power $g^k$ belonging to $F_2\times F_2$. Assuming this, we write $L\inc g\m Lg \inc\dots\inc g^{-k}Lg^{k}$.  By Lemma 1, we have $L=g^{-k}Lg^{k}$. Thus  $L= g\m Lg  =\langle L,h\rangle $, so $h\in L$ as required. 

Let $\Delta$ be the diagonal subgroup of $F_2\times F_2$. We let it act on the vector space  $\Q^4$ via the embedding of $F_2\times F_2$ into $GL(4,\Z)$. The only $\Delta$-invariant 2-planes are those of the form $H_\lambda=\{(x,\lambda x)\mid x\in\Q^2\}$, including the horizontal plane $H_0=\Q^2\times \{0\}$ and the vertical plane $H_\infty=\{0\}\times \Q^2$. For $\lambda\ne 0,\infty$, the setwise stabilizer of $H_\lambda$ for the action of  $F_2\times F_2$ on $\Q^4$ is $\Delta$.

Now $L$  properly contains $\Delta$, and   therefore  $H_0$, $H_\infty$ are the only $L$-invariant planes. Since $g$ conjugates $L$ to its overgroup $\langle L,h\rangle $, it must preserve the decomposition $\Q^4=\Q^2\oplus\Q^2$. Thus  $g^2\in GL(2,\Z) \times GL(2,\Z)$. Since $F_2\times F_2$ has finite index in $GL(2,\Z) \times GL(2,\Z)$ we deduce that some power of $g$ is in $F_2\times F_2$, as required.

\def\cprime{$'$}

\begin{flushleft}
 
 Laboratoire de Math\'ematiques Nicolas Oresme, %\\
Universit\'e de Caen et CNRS (UMR 6139), %)\\
BP 5186, %\\
F-14032 Caen Cedex, %\\
France\\

\emph{e-mail:} \texttt{levitt@math.unicaen.fr}\\
\end{flushleft}

\end{document}